\documentclass{amsart}
\usepackage{amssymb,amsmath,latexsym,epsfig}

\newtheorem{thm}{Theorem}[section]
\newtheorem{dfn}[thm]{Definition}
\newtheorem{cor}[thm]{Corollary}

\newtheorem{lemma}[thm]{Lemma}

\newcommand{\cl}{\hbox{\rm cl}}

\title[Lattice path matroids]{Lattice path matroids: enumerative aspects and Tutte polynomials}
\date{\today}
\author[Joseph E.~Bonin]
       {Joseph E.~Bonin}
\address[Joseph E.~Bonin]
{Department of Mathematics\\ The George Washington University\\
Washington, D.C. 20052, USA} \email [Joseph E.~Bonin]
{jbonin@gwu.edu}
\author[Anna de Mier]
{Anna de Mier}
\address[Anna de Mier]
{Departament de Matem\`atica Aplicada II\\ Universitat
Polit\`ecnica de Catalunya\\ Pau Gargallo 5, 08028\\ Barcelona,
Spain} \email [Anna de Mier] {demier@ma2.upc.es}
\author[Marc Noy]
{Marc Noy}
\address[Marc Noy]
{Departament de Matem\`atica Aplicada II\\ Universitat
Polit\`ecnica de Cata\-lunya\\ Pau Gargallo 5, 08028\\ Barcelona,
Spain} \email [Marc Noy] {noy@ma2.upc.es} 

\begin{document}

\begin{abstract}
Fix two lattice paths $P$ and $Q$ from $(0,0)$ to $(m,r)$ that use
East and North steps with $P$ never going above $Q$. We show that
the lattice paths that go from $(0,0)$ to $(m,r)$ and that remain
in the region bounded by $P$ and $Q$ can be identified with the
bases of a particular type of transversal matroid, which we call a
lattice path matroid. We consider a variety of enumerative aspects
of these matroids and we study three important matroid invariants,
namely the Tutte polynomial and, for special types of lattice path
matroids, the characteristic polynomial and the $\beta$ invariant.
In particular, we show that the Tutte polynomial is the generating
function for two basic lattice path statistics and we show that
certain sequences of lattice path matroids give rise to sequences
of Tutte polynomials for which there are relatively simple
generating functions. We show that Tutte polynomials of lattice
path matroids can be computed in polynomial time. Also, we obtain
a new result about lattice paths from an analysis of the $\beta$
invariant of certain lattice path matroids.
\end{abstract}

\maketitle

\section{Introduction}\label{sec:intro}

This paper develops a new connection between matroid theory and
enumerative combinatorics: with every pair of lattice paths $P$
and $Q$ that have common endpoints we associate a matroid in such
a way that the bases of the matroid correspond to the paths that
remain in the region bounded by $P$ and $Q$. These matroids, which
we call lattice path matroids, appear to have a wealth of
interesting and striking properties. In this paper we focus on the
enumerative aspects of lattice path matroids, including the study
of important matroid invariants like the Tutte and the
characteristic polynomials. Structural aspects of lattice path
matroids and their relation with other families of matroids will
be the subject of a forthcoming paper~\cite{lpm2}.

Lattice path matroids provide a bridge between matroid theory and
the theory of lattice paths that, as we demonstrate here and
in~\cite{lpm2}, can lead to a mutually enriching relationship
between the two subjects. One example starts with the path
interpretation we give for each coefficient of the Tutte
polynomial of a lattice path matroid.  Computing the Tutte
polynomial of an arbitrary matroid is known to be $\#$P-complete;
the same is true even within special classes such as graphic or
transversal matroids. However, by using the path interpretation of
the coefficients in the case of lattice path matroids, we show
that the Tutte polynomial of such a matroid can be computed in
polynomial time. On the lattice path side, as we illustrate in
Section~\ref{sec:beta}, this interpretation of the coefficients
along with easily computed examples of the Tutte polynomial can
suggest new theorems about lattice paths.

Relatively little matroid theory is required to understand this
paper and what is needed is sketched in the first part of
Section~\ref{sec:back}. We follow the conventional notation for
matroid theory as found in~\cite{ox}. A few topics of matroid
theory of a more specialized nature (the Tutte polynomial, the
broken circuit complex, the characteristic polynomial, and the
$\beta$ invariant) are presented in the sections in which they
play a role. The last part of Section~\ref{sec:back} outlines the
basic facts on lattice path enumeration that we use.

Lattice path matroids, the main topic of this paper, are defined
in Section~\ref{sec:lpm}, where we also identify their bases with
lattice paths (Theorem~\ref{thm:bases}). We introduce special
classes of lattice path matroids, among which are the Catalan
matroids and, more generally, the $k$-Catalan matroids, for which
the numbers of bases are the Catalan numbers and the $k$-Catalan
numbers. We also treat some basic enumerative results for lattice
path matroids and prove several structural properties of these
matroids that play a role in enumerative problems that are
addressed later in the paper. Counting connected lattice path
matroids on a given number of elements is the topic of
Section~\ref{sec:enum}.

The next four sections consider matroid invariants in the case of
lattice path matroids. Section~\ref{sec:tutte} gives a lattice
path interpretation of each coefficient of the Tutte polynomial of
a lattice path matroid (Theorem~\ref{active}) as well as
generating functions for the Tutte polynomials of the sequence of
$k$-Catalan matroids (Theorem~\ref{thm:tuttegf}) and, from that, a
formula for each coefficient of each of these Tutte polynomials
(Theorem~\ref{thm:tuttecoefs}).  In Section~\ref{sec:computing} we
give an algorithm for computing the Tutte polynomial of any
lattice path matroid in polynomial time; we provide a second
method of computation that applies for certain classes of lattice
path matroid and which, although limited in scope, is particularly
simple to implement on a computer.  In Section~\ref{sec:char} we
show that the broken circuit complex of a lattice path matroid is
the independence complex of another lattice path matroid and we
develop formulas for the coefficients of the characteristic
polynomial for special classes of lattice path matroids.
Section~\ref{sec:beta} shows that $k$ times the Catalan number
$C_{kn-1}$ counts lattice paths of a special type
(Theorem~\ref{thm:betaud}); the key to discovering this result was
looking at a particular coefficient (the $\beta$ invariant) of the
Tutte polynomials of certain lattice path matroids.

The final section connects lattice path matroids with a problem of
current interest in enumerative combinatorics, namely, the
$(k+l,l)$-tennis ball problem.

We use the following common notation: $[n]$ denotes the set
$\{1,2,\dots,n\}$ and $[m,n]$ denotes the set $\{m,n+1,\dots,
n\}$. We follow the  convention in matroid theory of writing
$X\cup e$ and $X-e$ in place of $X\cup \{e\}$ and $X-\{e\}$.

\section{Background}\label{sec:back}

In this section we introduce the concepts of matroid theory that
are needed in this paper. For a thorough introduction to the
subject we refer the reader to Oxley~\cite{ox}; the proofs we omit
in this section can be found there, mostly in Chapters $1$ and
$2$.  We conclude this section with the necessary background on
the enumerative theory of lattice paths.

\begin{dfn}
A \emph{matroid} is a pair $\bigl(E(M), \mathcal{B}(M)\bigr)$
consisting of a finite set $E(M)$ and a collection
$\mathcal{B}(M)$ of subsets of $E(M)$ that satisfy the following
conditions:
\begin{itemize}
\item[(B1)] $\mathcal{B}(M)\ne\emptyset$,
\item[(B2)] $\mathcal{B}(M)$ is an antichain, that is, no set in
$\mathcal{B}(M)$ properly contains another set in
$\mathcal{B}(M)$, and
\item[(B3)] for each pair of distinct sets $B,B'$ in $\mathcal{B}(M)$
and for each element $x\in B-B'$, there is an element $y\in B'-B$
such that $(B-x)\cup y$ is in $\mathcal{B}(M)$.
\end{itemize}
\end{dfn}

The set $E(M)$ is  the \emph{ground set} of $M$ and the sets in
$\mathcal{B}(M)$ are  the \emph{bases} of $M$. Subsets of bases
are {\em independent sets}; the collection of independent sets of
$M$ is denoted $\mathcal{I}(M)$. Sets that are not independent are
\emph{dependent}. A \emph{circuit} is a minimal dependent set. If
$\{x\}$ is a circuit, then $x$ is  a \emph{loop}. Thus, no basis
of $M$ can contain a loop. An element that is contained in every
basis is an \emph{isthmus}.

It is easy to show that all bases of $M$ have the same
cardinality.  More generally, for any subset $A$ of $E(M)$ all
maximal independent subsets of $A$ have the same cardinality;
$r(A)$, the {\em rank} of $A$, denotes this common cardinality. If
several matroids are under consideration, we may use $r_M(A)$ to
avoid ambiguity. In place of $r(E(M))$, we write $r(M)$.

The \emph{closure} of a set $A\subseteq E(M)$ is defined as
$$\cl(A)=\{x\in E(M): r(A\cup x)=r(A)\}.$$ A set $F$ is a
\emph{flat} if $\cl (F)=F$.  The flats of a matroid, ordered by
inclusion, form a geometric lattice.

It is well-known that matroids can be characterized in terms of
each of the following objects: the independent sets, the dependent
sets, the circuits, the rank function, the closure operator, and
the flats (see Sections~1.1--1.4 of~\cite{ox}).

\smallskip

\noindent\textsc{Example.} A matroid of rank $r$ is a
\emph{uniform matroid} if  all $r$-element subsets of the ground
set are bases. There is, up to isomorphism, exactly one uniform
matroid of rank $r$ on an $m$-element set; this matroid is denoted
$U_{r,m}$.

\smallskip

One fundamental concept in matroid theory is duality. Given a
matroid $M$, its \emph{dual matroid} $M^*$ is the matroid on
$E(M)$ whose set of bases is given by $$\mathcal{B}(M^*)=\{E(M)-B:
B\in \mathcal{B}(M)\}.$$  A matroid is \emph{self-dual} if it is
isomorphic to its dual; a matroid is \emph{identically self-dual}
if it is equal to its dual. For example, the dual of the uniform
matroid $U_{r,m}$ is the uniform matroid $U_{m-r,m}$. The matroid
$U_{r,2r}$ is identically self-dual. For any matroid $M$, the
element $x$ is a loop of $M$ if and only if $x$ is an isthmus of
the dual $M^*$.

This paper investigates a special class of transversal matroids.
Let $\mathcal{A}=(A_j: j\in J)$ be a set system, that is, a
multiset of subsets of a finite set $S$. A \emph{transversal} (or
system of distinct representatives) of $\mathcal{A}$ is a set
$\{x_j:j\in J\}$ of $|J|$ distinct elements such that $x_j\in A_j$
for all $j$ in $J$. A \emph{partial transversal} of $\mathcal{A}$
is a transversal of a set system of the form $(A_k:k\in K)$ with
$K$ a subset of $J$. The following theorem is a fundamental result
due to Edmonds and Fulkerson.

\begin{thm}
The  partial transversals of a set system $\mathcal{A}=(A_j: j\in
J)$ are the independent sets of a matroid on $S$.
\end{thm}

A \emph{transversal matroid} is a matroid whose independent sets
are the partial transversals of some set system $\mathcal{A}=(A_j:
j\in J)$; we say that $\mathcal{A}$ is a \emph{presentation} of
the transversal matroid. The bases of a transversal matroids are
the maximal partial transversals of $\mathcal{A}$. For more on
transversal matroids see~\cite[Section 1.6]{ox}.

Given two matroids $M_1,M_2$ on disjoint ground sets, their
\emph{direct sum} is the matroid $M_1\oplus M_2$ with ground set
$E(M_1)\cup E(M_2)$ whose collection of bases is
$$\mathcal{B}(M_1\oplus M_2)= \{B_1\cup B_2: B_1\in
\mathcal{B}(M_1) \text{ and } B_2\in \mathcal{B}(M_2)\}.$$ It is
easy to check that the lattice of flats of $M_1\oplus M_2$ is
isomorphic to the direct product (or cartesian product) of the
lattice of flats of $M_1$ and that of $M_2$. A matroid $M$ is
\emph{connected} if it is not a direct sum of two nonempty
matroids. Note that connected matroids with at least two elements
have neither loops nor isthmuses.

We say that the matroid $M\oplus U_{1,1}$ is formed by
\emph{adding an isthmus} to $M$. In the case that the ground set
of the uniform matroid $U_{1,1}$ is $e$, we shorten this notation
to $M\oplus e$ if there is no danger of ambiguity. Of course, $e$
is an isthmus of $M\oplus e$.

There is a well developed theory of extending matroids by single
elements~\cite[Section 7.2]{ox}.  The case that is relevant to
this paper is that of free extension, which consists of adding an
element to the matroid as independently as possible without
increasing the rank. Precisely, the {\em free extension} $M+e$ of
$M$ is the matroid on $E(M)\cup e$ whose collection of independent
sets is given as follows: $$\mathcal{I}(M+ e)=\mathcal{I}(M)\cup
\{I\cup e: I\in \mathcal{I}(M) \text{ and } |I|<r(M)\}.$$ The
bases of $M+e$, where $M$ has rank $r$, are the bases of $M$
together with the sets of the form $I\cup e$, where $I$ is an
$(r-1)$-element independent set of $M$. Equivalently, the rank
function of $M+e$ is given by the following equations: for $X$ a
subset of $E(M)$, $$r_{M+e}(X) = r_M(X)$$ and
\begin{equation*}
r_{M+e}(X\cup e) =    \left\{
     \begin{array}{ll}
      r_M(X)+1, &\mbox{if $r_M(X)<r(M)$;}\\
      r(M), &\mbox{otherwise.}
    \end{array}\notag
   \right.
\end{equation*}

The particular matroids of interest in this paper arise from
lattice paths, to which we now turn. We consider two kinds of
lattice paths, both of which are in the plane. Most of the lattice
paths we consider use steps $E=(1,0)$ and $N=(0,1)$; in several
cases it is more convenient to use lattice paths with steps
$U=(1,1)$ and $D=(1,-1)$.  The letters are abbreviations of East,
North, Up, and Down. We will often treat lattice paths as words in
the alphabets $\{E,N\}$ or $\{U,D\}$, and we will use the notation
$\alpha^n$ to denote the concatenation of $n$ letters, or strings
of letters, $\alpha$. If $P=s_1s_2 \ldots s_n$ is a lattice path,
then its \emph{reversal} is defined as $P^{\rho} = s_n s_{n-1}
\dots s_1$. The {\em length} of a lattice path $P=s_1s_2 \ldots
s_n$ is $n$, the number of steps in $P$.

Here we recall the facts we need about the enumeration of lattice
paths; the proofs of the following lemmas can be found in Sections
3 to 5 of the first chapter of~\cite{moh}. The most basic
enumerative results about lattice paths are those in the following
lemma.

\begin{lemma}\label{countpath}
For a fixed positive integer $k$, the number of lattice paths from
$(0,0)$ to $(kn,n)$ that use steps $E$ and $N$ and that never pass
above the line $y=x/k$ is the $n$-th $k$-Catalan number $$C^k_n =
\frac{1}{kn+1}\binom{(k+1)n}{n}.$$ In particular, the number of
paths from $(0,0)$ to $(n,n)$ that never pass above the line $y=x$
is the $n$-th Catalan number $$C_n = \frac{1}{n+1}\binom{2n}{n}.$$
\end{lemma}

We also use the following result, which generalizes
Lemma~\ref{countpath}. For $k=1$ the numbers displayed in
Lemma~\ref{ballot} are called the \emph{ballot numbers}.

\begin{lemma}\label{ballot}
For $m\geq kn\geq 0$, the number of lattice paths from $(0,0)$ to
$(m,n)$ with steps $E$ and $N$ that never go above the line
$y=x/k$ is $$\frac{m-kn+1}{m+n+1}\binom{m+n+1}{n}.$$
\end{lemma}

The next lemma treats paths in the alphabet $\{U,D\}$; the first
assertion, which concerns what are usually called Dyck paths, is
equivalent to the second part of Lemma~\ref{countpath} by the
obvious identification of the alphabets.

\begin{lemma}\label{lemma:dyck}
(i) The number of paths from $(0,0)$ to $(2n,0)$ with steps $U$
and $D$ that never pass below the $x$-axis is the $n$-th Catalan
number $C_n$.
\\
(ii) The number of paths of $n$ steps in the alphabet $\{U,D\}$
that start at $(0,0)$ and never pass below the $x$-axis (not
necessarily ending on the $x$-axis) is $\binom{n}{\lceil n/2
\rceil}$.
\end{lemma}

The following result will be used to count certain types of
lattice paths.

\begin{lemma}\label{lem:coeffcat}
Let $$C(z) = \sum_{n\geq 0} \frac{1}{kn+1}\binom{(k+1)n}{n} z^n$$
be the generating function for the $k$-Catalan numbers. The
coefficient of $z^t$ in $C(z)^j$ is
$$\frac{j}{t}\binom{(k+1)t+j-1}{t-1}.$$
\end{lemma}

\section{Lattice path matroids}\label{sec:lpm}

In this section we define lattice path matroids as well as several
important subclasses.  Later sections of this paper develop much
of the enumerative theory for lattice path matroids in general and
this theory is pushed much further for certain special families of
lattice path matroids.

\begin{dfn}\label{def:lpm}
Let $P=p_1p_2\dots p_{m+r}$ and $Q=q_1q_2\dots q_{m+r}$ be two
lattice paths from $(0,0)$ to $(m,r)$ with $P$ never going above
$Q$. Let $\{p_{u_1},\dots,p_{u_r}\}$ be the set of North steps of
$P$ with $u_1<u_2<\cdots<u_r$; similarly, let
$\{q_{l_1},\dots,q_{l_r}\}$ be the set of North steps of $Q$ with
$l_1<l_2<\cdots<l_r$. Let $N_i$ be the interval $[l_i,u_i]$ of
integers. Let $M[P,Q]$ be the transversal matroid that has ground
set $[m+r]$ and presentation $(N_i: i\in [r])$; the pair $(P,Q)$
is a \emph{presentation of $M[P,Q]$}. A \emph{lattice path
matroid} is a matroid $M$ that is isomorphic to $M[P,Q]$ for some
such pair of lattice paths $P$ and $Q$.
\end{dfn}

Several examples of lattice path matroids are given after
Theorem~\ref{thm:bases}, which identifies the bases of these
matroids in terms of lattice paths. To avoid needless repetition,
throughout the rest of the paper we assume that the lattice paths
$P$ and $Q$ are as in Definition~\ref{def:lpm}.

We think of $1,2,\dots, m+r$ as the first step, the second step,
etc. Observe that the set $N_i$ contains the steps that can be the
$i$-th North step in a lattice path from $(0,0)$ to $(m,r)$ that
remains in the region bounded by $P$ and $Q$. When thought of as
arising from the particular presentation using bounding paths $P$
and $Q$, the elements of the matroid are ordered in their natural
order, i.e., $1<2<\cdots<m+r$; we will frequently use this order
throughout the paper. However, this order is not inherent in the
matroid structure; the elements of a lattice path matroid
typically can be linearly ordered in many ways so as to correspond
to steps in lattice paths.  (This point will be addressed in
greater detail in~\cite{lpm2}.)

We associate a lattice path $P(X)$ with each subset $X$ of the
ground set of a lattice path matroid as specified in the next
definition.

\begin{dfn}\label{path}
Let $X$ be a subset of the ground set $[m+r]$ of the lattice path
matroid $M[P,Q]$. The lattice path $P(X)$ is the word $$s_1 s_2
\dots s_{m+r}$$ in the alphabet $\{E,N\}$ where
\begin{eqnarray}
 s_i =
   \left\{
     \begin{array}{ll}
      N, &\mbox{if $i\in X$;}\\
      E, &\mbox{otherwise.}
    \end{array}\notag
   \right.
\end{eqnarray}
\end{dfn}

Thus, the path $P(X)$ is formed by taking the elements of $M[P,Q]$
in the natural linear order and replacing each by a North step if
the element is in $X$ and by an East step if the element is not in
$X$.

The fundamental connection between the transversal matroid
$M[P,Q]$ and the lattice paths that stay in the region bounded by
$P$ and $Q$ is the following theorem which says that the bases of
$M[P,Q]$ can be identified with such lattice paths.

\begin{thm}\label{thm:bases}
A subset $B$ of $[m+r]$ with $|B|=r$ is a basis of $M[P,Q]$ if and
only if the associated lattice path $P(B)$ stays in the region
bounded by $P$ and $Q$.
\end{thm}

\begin{proof}
Let $B$ be $\{b_1,\dots, b_r\}$ with $b_1<b_2<\cdots<b_r$ in the
natural order. Suppose first that $B$ is a basis of $M[P,Q]$, that
is, a transversal of $(N_i:i\in [r])$. The conclusion will follow
if we prove that $b_i$ is in $N_i$. Assume, to the contrary,
$b_i\not\in N_i$. Since either $b_i<l_i$ or $b_i>u_i$, we obtain
the following contradictions: in the first case, the set
$\{b_1,b_2\dots,b_i\}$ must be a transversal of
$(N_1,N_2,\dots,N_{i-1})$; in the second,
$\{b_i,b_{i+1},\dots,b_r\}$ must be a transversal of
$(N_{i+1},N_{i+2},\dots,N_r)$.

Conversely, if the lattice path $P(B)$ goes neither below $P$ nor
above $Q$, then for every $i$ we have that $b_i$, the $i$-th North
step of $P(B)$, satisfies $l_i\leq b_i\leq u_i$, and hence that
$B$ is a transversal of $(N_i:i\in [r])$.
\end{proof}

\begin{cor}\label{cor:bases}
The number of basis of $M[P,Q]$ is the number of lattice paths
from $(0,0)$ to $(m,r)$ that go neither below $P$ nor above $Q$.
\end{cor}

\begin{figure}
\begin{center}
\epsfxsize 5.5truein \epsffile{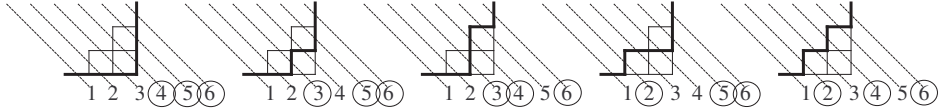}
\end{center}
\caption{The bases $\{4,5,6\}$, $\{3,5,6\}$, $\{3,4,6\}$,
$\{2,5,6\}$, and $\{2,4,6\}$ of a lattice path matroid represented
as the North steps of lattice paths.}\label{catex}
\end{figure}

Figure~\ref{catex} illustrates Theorem~\ref{thm:bases}.  In this
example we have $N_1=\{2,3,4\}$,  $N_2=\{4,5\}$, and $N_3=\{6\}$.
There are five bases of this transversal matroid.  Note that $1$
is a loop and $6$ is an isthmus.

\smallskip

\noindent\textsc{Example}. For the lattice paths $P=E^mN^r$ and
$Q=N^rE^m$, every $r$-subset of $[m+r]$ is a basis of $M[P,Q]$.
Thus, the uniform matroid $U_{r,m+r}$ is a lattice path matroid.

\smallskip

\begin{figure}
\begin{center}
\epsfxsize 4.0truein \epsffile{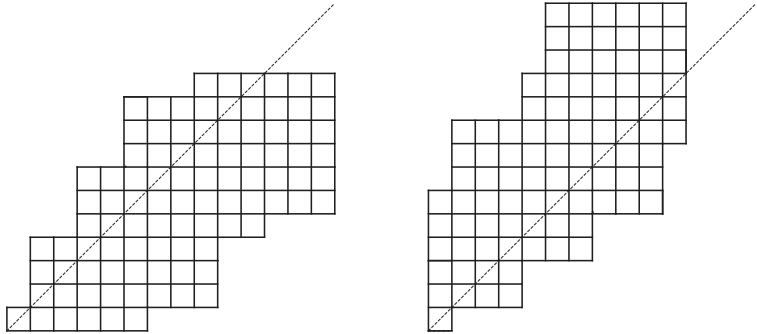}
\end{center}
\caption{Presentations of a lattice path matroid and its
dual.}\label{dual}
\end{figure}

Recall that the bases of the dual $M^*$ of a matroid $M$ are the
set complements of the bases of $M$ with respect to the ground set
$E(M)$. Thus, for a lattice path matroid $M$, the bases of the
dual matroid correspond to the East steps in lattice paths.
Reflecting a lattice path presentation of $M$ about the line $y=x$
shows that the dual $M^*$ is also a lattice path matroid. (See
Figure~\ref{dual}.) This justifies the following theorem.

\begin{thm}\label{thm:dual}
The class of lattice path matroids is closed under matroid
duality.
\end{thm}

Note that a $180^\circ$ rotation of the region bounded by $P$ and
$Q$, translated to start at $(0,0)$, yields the same matroid
although the labels on the elements are reversed. Thus the lattice
path matroids  $M[P,Q]$ and $M[Q^{\rho},P^{\rho}]$ are isomorphic.
It follows, for example, that the lattice path matroid in
Figure~\ref{catex} is self-dual; note that this matroid is not
identically self-dual since the loop $1$ and the isthmus $6$ in
the matroid are, respectively, an isthmus and a loop in the dual.

\begin{figure}
\begin{center}
\epsfxsize 4.0truein \epsffile{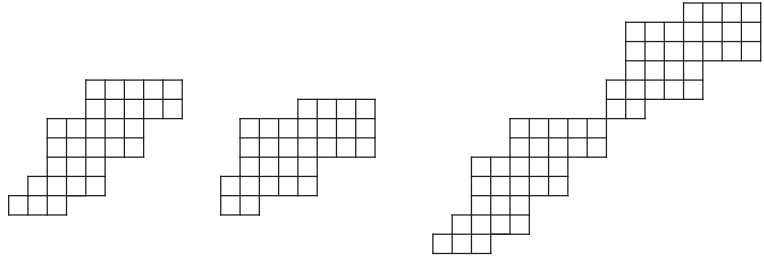}
\end{center}
\caption{Presentations of two lattice path matroids and their
direct sum.}\label{sum}
\end{figure}

Figure~\ref{sum} illustrates the next result.  The proof is
immediate  from Theorem~\ref{thm:bases} and the definition of
direct sums.

\begin{thm}\label{dsum}
The class of lattice path matroids is closed under direct sums.
Furthermore, the lattice path matroid $M[P,Q]$ is connected if and
only if the bounding lattice paths $P$ and $Q$ meet only at
$(0,0)$ and $(m,r)$.
\end{thm}

We now turn to a special class of lattice path matroids, the
generalized Catalan matroids, as well as to various subclasses
that exhibit a structure that is simpler than that of typical
lattice path matroids. Later sections of this paper will give
special attention to these classes since the simpler structure
allows us to obtain more detailed enumerative results.

\begin{dfn}
A lattice path matroid $M$ is a \emph{generalized Catalan matroid}
if there is a presentation $(P,Q)$ of $M$ with $P=E^m N^r$. In
this case we simplify the notation $M[P,Q]$ to $M[Q]$. If in
addition the upper path $Q$ is $(E^k N^l)^n$ for some positive
integers $k,l,$ and $n$, we say that $M$ is the
\emph{$(k,l)$-Catalan matroid} $M_n^{k,l}$.  In place of
$M_n^{k,1}$ we write $M_n^k$; such matroids are called
\emph{$k$-Catalan matroids}.  In turn, we simplify the notation
$M_n^1$ to $M_n$; such matroids are called \emph{Catalan
matroids}.
\end{dfn}

\begin{figure}
\begin{center}
\epsfxsize 4.0truein \epsffile{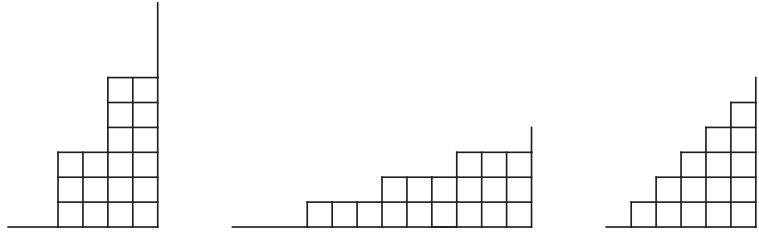}
\end{center}
\caption{Presentations of the rank nine matroid $M^{2,3}_3$, the
$3$-Catalan matroid $M^3_4$ of rank four, and the rank six Catalan
matroid $M_6$. }\label{gcmex}
\end{figure}

Figure~\ref{gcmex} gives presentations of a $(2,3)$-Catalan
matroid, a $3$-Catalan matroid, and a Catalan matroid. These
matroids have, respectively, two loops and three isthmuses, three
loops and one isthmus, and a single loop and isthmus.

Note that $(k,l)$-Catalan matroids have isthmuses and loops;
specifically, the elements $1,\dots,k$ are the loops and
$(k+l)n-l+1,(k+l)n-l+2,\dots,(k+l)n$ are the isthmuses of
$M^{k,l}_n$. Also, observe that for the $k$-Catalan matroid
$M^k_n$, Theorem~\ref{thm:bases} can be restated by saying that an
$n$-element subset $B$ of $[(k+1)n]$ is a basis of $M^k_n$ if and
only if its associated lattice path $P(B)$ does not go above the
line $y=x/k$.

We next note an immediate consequence of Corollary~\ref{cor:bases}
and Lemma~\ref{countpath}. As we will see in
Section~\ref{sec:tennis}, there is no known formula that leads to
a similar result for $(k,l)$-Catalan matroids.

\begin{cor}
The number of bases of the  $k$-Catalan matroid $M_n^k$ is the
 $k$-Catalan number $C^k_n$. In particular, the number of
bases of the Catalan matroid $M_n$ is the  Catalan number $C_n$.
\end{cor}

The comments before and after Theorem~\ref{thm:dual}, including
that about $180^\circ$ rotations of presentations, give the
following result.

\begin{thm}\label{thm:dualkl}
The dual of the  $(k,l)$-Catalan matroid $M^{k,l}_n$ is the
$(l,k)$-Catalan matroid $M^{l,k}_n$.  Thus, the matroid
$M^{k,k}_n$, and in particular the Catalan matroid $M_n$, is
self-dual but not identically self-dual.
\end{thm}

We turn to lattice path descriptions of circuits and independent
sets in generalized Catalan matroids. Recall from
Definition~\ref{path} that we associate a lattice path $P(X)$ with
each subset $X$ of the ground set $[m+r]$ of the lattice path
matroid $M[P,Q]$.  Of course, only sets of $r$ elements give paths
that end at $(m,r)$.

\begin{thm}\label{thm:circuits}
A subset $C$ of $[m+r]$ is a circuit of the generalized Catalan
matroid $M[Q]$ if and only if for the largest element $i$ of $C$,
the $i$-th step of $P(C)$ is the only North step of $P(C)$ above
$Q$.
\end{thm}

\begin{proof}
First assume that for the largest element $i$ of $C$, the $i$-th
step of $P(C)$ is the only North step of $P(C)$ above $Q$.  It is
clear that for any superset $X$ of $C$, the $i$-th step of $P(X)$
is also above $Q$.  Thus, $C$ is not contained in any basis and so
is dependent.  Note that for any element $c$ in $C$, the lattice
path $P(C-c)$ has no steps above $Q$; also, the path that follows
$P(C-c)$ to the line $x=m$ and then goes directly North to $(m,r)$
is a lattice path that never goes above $Q$ and so corresponds to
a basis that contains $C-c$, specifically, the basis $(C-c)\cup Y$
where $Y$ contains the last $r-(|C|-1)$ elements in $[m+r]$. Thus,
every proper subset of $C$ is independent. Therefore, $C$ is a
circuit.

Conversely, assume that $C$ is a circuit. By the same type of
argument as in the second half of the last paragraph, it is clear
that $P(C)$ must have at least one North step that goes above $Q$;
since $C$ is a minimal dependent set, it is clear that this step
must correspond to the greatest element of $C$.
\end{proof}

\begin{cor}
The number of $i$-element circuits in the Catalan matroid $M_n$ is
the Catalan number $C_{i-1}$.
\end{cor}

\begin{proof}
From the last theorem, it follows that the lattice path $P(C)$
associated with an $i$-element circuit $C$ can be decomposed as
follows: a lattice path from $(0,0)$ to $(i-1,i-1)$ that does not
go above the line $y=x$, followed by one North step above the line
$y=x$, followed by only East steps.  Conversely, any such path
corresponds to an $i$-element circuit. From this the result
follows.
\end{proof}

From Theorem~\ref{thm:circuits} we also get the following result.

\begin{cor}\label{cor:indpaths}
The independent sets in the generalized Catalan matroid $M[Q]$ are
precisely the subsets $X$ of $[m+r]$ such that the associated
lattice path $P(X)$ never goes above the bounding lattice path
$Q$.
\end{cor}

From this result it follows that for $k$-Catalan matroids, the
paths that correspond to independent sets of a given size are
precisely those given by Lemma~\ref{ballot}.

\begin{cor}
The number of independent sets of size $i$ in the $k$-Catalan
matroid $M^k_n$ is
$$\frac{(k+1)(n-i)+1}{(k+1)n+1}\binom{(k+1)n+1}{i}.$$
\end{cor}

Generalized Catalan matroids have previously appeared in the
matroid theory literature~\cite{number}; they were also studied,
from a very different perspective, in~\cite{opr} and they were
recently rediscovered in yet another context in~\cite{ardila}. It
can be shown that generalized Catalan matroids are exactly the
minors of Catalan matroids~\cite{lpm2}. We conclude this section
with yet another perspective by showing that generalized Catalan
matroids are precisely the matroids that can be constructed from
the empty matroid by repeatedly adding isthmuses and taking free
extensions. Theorem~\ref{thm:catrec} can be generalized to all
lattice path matroids; the generalization uses more matroid
theory, in particular, more general types of extensions than free
extensions, so this result will appear in~\cite{lpm2}.  We present
the special case here since in the last part of
Section~\ref{sec:computing} we will use this result to give simple
and efficient algebraic rules to compute the Tutte polynomial of
any generalized Catalan matroid.

\begin{thm}\label{thm:catrec}
Let $Q=q_1 q_2\dots q_{m+r}$ be a word of length $m+r$ in the
alphabet $\{E,N\}$. Let $M^0$ be the empty matroid and
define$$M^i=\left\{
\begin{array}{ll} M^{i-1} + i, & \mathrm{if \ } q_i=E; \\ M^{i-1}
\oplus i, & \mathrm{if \ } q_i=N. \end{array} \right.$$ Then
$M^{m+r}$ and the generalized Catalan matroid $M[Q]$ are equal.
\end{thm}

\begin{proof}
Let $Q_i$ be the initial segment $q_1 q_2 \dots q_{i}$ of $Q$, let
$R_i$ be the region determined by the bounding paths of $M[Q_i]$,
and let the paths that correspond to bases of $M[Q_i]$ go from
$(0,0)$ to $(m_i,r_i)$. We prove the equality $M^i = M[Q_i]$ by
induction on $i$. Both $M^0$ and $M[Q_0]$ are the empty matroid.
Assume $M^{i-1} = M[Q_{i-1}]$. Assume first that $q_i$ is $N$, so
$i$ is an isthmus of $M^i$.  Thus we need to show that the bases
of $M[Q_i]$ are precisely the sets of the form $B\cup i$ where $B$
is a basis of $M[Q_{i-1}]$, which is clear from
Theorem~\ref{thm:bases} since the bounding paths for $M[Q_i]$ have
a common last ($i$-th) North step.  Now suppose that $q_i$ is $E$.
Note the equality $r_i = r_{i-1}.$ Lattice paths in the region
$R_i$ from $(0,0)$ to $(m_i,r_i)$ are of two types: those in which
the final step is North, and so correspond to sets of the form
$I\cup i$ where $I$ is an independent set of size $r_{i-1}-1$ in
$M[Q_{i-1}]$; those in which the final step is East, and so
correspond to bases of $M[Q_{i-1}]$.  From this and the basis
formulation of free extensions, the equality $M^i = M[Q_i]$
follows.
\end{proof}

\section{Enumeration of lattice path matroids}\label{sec:enum}

In this section we give a formula for the number of connected
lattice matroids on a given number of elements up to isomorphism;
to make the final result slightly more compact, we let the number
of elements be $n+1$. The proof has two main ingredients, the
first of which is the following result from~\cite{lpm2}. (Recall
that $P^{\rho}$ denotes the reversal $s_{n+1} s_n \dots s_1$ of a
lattice path $P=s_1s_2 \ldots s_n s_{n+1}$.)

\begin{lemma}\label{lemma:iso}
Two connected lattice path matroids $M[P,Q]$ and $M[P',Q']$ are
isomorphic if and only if either $P'=P$ and $Q'=Q$, or
$P'=Q^{\rho}$ and $Q'=P^{\rho}$.
\end{lemma}

\begin{figure}
\begin{center}
\epsfxsize 3.5truein \epsffile{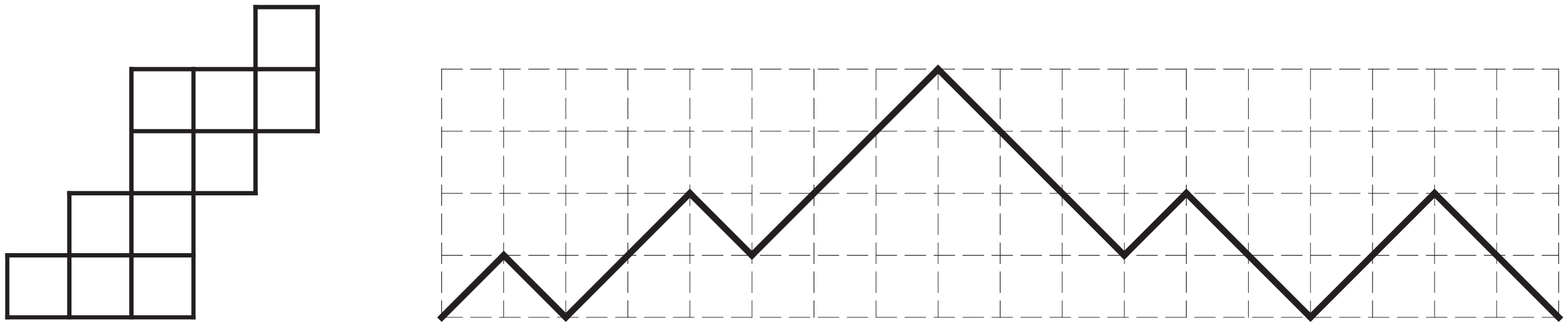}
\end{center}
\caption{A parallelogram polyomino and its associated Dyck path.
}\label{fig:dyck}
\end{figure}

The second main ingredient is the following bijection, going back
at least to P\'olya, between the pairs of lattice paths of length
$n+1$ that intersect only at their endpoints and the Dyck paths of
length $2n$. (See, for example,~\cite{dv}.) A pair $(P,Q)$ of
nonintersecting lattice paths from $(0,0)$ to $(m,r)$ can be
viewed as the special type of polyomino that in~\cite{dv} is
called a parallelogram polyomino.  Associate two sequences
$(a_1,\dots,a_m)$ and $(b_1, \dots, b_{m-1})$ of integers with
such a polyomino: $a_i$ is the number of cells of the $i$-th
column of the polyomino (columns are scanned from left to right)
and $b_i+1$ is the number of cells of column $i$ that are adjacent
to cells of column $i+1$. Since the paths are nonintersecting,
each $b_i$ is nonnegative. Now associate to $(P,Q)$ the Dyck path
$\pi$ having $m$ peaks at heights $a_1,\dots,a_m$ and $m-1$
valleys at heights $b_1,\dots,b_{m-1}$. Figure~\ref{fig:dyck}
shows a polyomino and its associated Dyck path; the corresponding
sequences for this polyomino are $(1,2,4,2,2)$ and $(0,1,1,0)$. It
can be checked that the correspondence $(P,Q) \mapsto \pi$ is
indeed a bijection. Hence the number of such pairs $(P,Q)$ of
lattice paths of length $n+1$ is the Catalan number $C_n$.

Note that $C_n$ is not the number of connected lattice path
matroids on $n+1$ elements since different pairs of paths can give
the same matroid. According to Lemma~\ref{lemma:iso}, this happens
only for a pair $(P,Q)$ and its reversal $(Q^{\rho},P^{\rho})$, so
we need to find the number of pairs $(P,Q)$ for which $(P,Q) =
(Q^{\rho},P^{\rho})$. It is immediate to check that $(P,Q) =
(Q^{\rho},P^{\rho})$ if and only if the corresponding Dyck path
$\pi$ is symmetric with respect to its center or, in other words,
is equal to its reversal. Since a symmetric Dyck path of length
$2n$ is determined by its first $n$ steps, the number of such
paths is given in part (ii) of Lemma~\ref{lemma:dyck}. From the
number $C_n$ we obtained in the last paragraph we must subtract
half the number of nonsymmetric Dyck paths, thus giving the
following result.

\begin{thm}
The number of connected lattice path matroids on $n+1$ elements up
to isomorphism is
 $$
    C_n - \frac{1}{2}\left(C_n - \binom{n}{\lceil n/2 \rceil}\right)
     = \frac{1}{2} C_n + \frac{1}{2} \binom{n}{\lceil n/2 \rceil}.
 $$
\end{thm}

This number is asymptotically of order $O(4^n)$. Since it is known
that the number of transversal matroids on $n$ elements grows like
$c^{n^2}$ for some constant $c$ (see~\cite{affine}), it follows
that the class of lattice path matroids is rather small with
respect to the class of all transversal matroids.

We remark that the total number of lattice path matroids
(connected or not) on $k$ elements is the number of multisets of
connected lattice path matroids, the sum of whose cardinalities is
$k$. A generating function for these numbers can be derived using
standard tools; however, the result does not seem to admit a
compact form so we omit it.

\section{Tutte polynomials}\label{sec:tutte}
The Tutte polynomial is one of the most widely studied matroid
invariants. From the Tutte polynomial one obtains, as special
evaluations, many other important polynomials, such as the
chromatic and flow polynomials of a graph, the weight enumerator
of a linear code, and the Jones polynomial of an alternating knot.
(See~\cite{bo,ckcc} for many of the numerous occurrences of this
polynomial in combinatorics, in other branches of mathematics, and
in other sciences.) In this section, after reviewing the
definition of the Tutte polynomial, we show that for lattice path
matroids this polynomial is the generating function for two basic
lattice path statistics.  We use this lattice path interpretation
of the Tutte polynomial to give a formula for the generating
function $\sum_{n\geq 0} t(M^k_n;x,y) z^n$ for the sequence of
Tutte polynomials $t(M^k_n;x,y)$ of the $k$-Catalan matroids.
Using this generating function, we then derive a formula for each
coefficient of the Tutte polynomial $t(M^k_n;x,y)$.

The Tutte polynomial $t(M;x,y)$ of a matroid $M$ is most briefly
defined as follows:
\begin{equation}\label{tdef}
t(M;x,y) = \sum_{A \subseteq E(M)}
(x-1)^{r(M)-r(A)}(y-1)^{|A|-r(A)}.
\end{equation}
However, for our work the formulation in terms of internal and
external activities, which we review below, will prove most
useful. For a proof of the equivalence of these definitions (and
that the formulation in terms of activities is well-defined), see,
for example~\cite{shell}.

Fix a linear order $<$ on $E(M)$ and let $B$ be a basis of $M$. An
element $e\not\in B$ is \emph{externally active with respect to
$B$} if there is no element $y$ in $B$ with $y<e$ for which
$(B-y)\cup e$ is a basis.  An element $b\in B$ is \emph{internally
active with respect to $B$} if there is no element $y$ in $E(M)-B$
with $y<b$ for which $(B-b)\cup y$ is a basis. The {\em internal
(external) activity} of a basis is the number of elements that are
internally (externally) active with respect to that basis. We
denote the activities of a basis $B$ by $i(B)$ and $e(B)$. Note
that $i(B)$ and $e(B)$ depend not only on $B$ but also on the
order $<$. The following lemma is well-known and easy to prove.

\begin{lemma}\label{dualact}
Let the elements of a matroid $M$ and its dual $M^*$ be ordered
with the same linear ordering.  An element $e$ is internally
active with respect to the basis $B$ of $M$ if and only if $e$ is
externally active with respect to the basis $E(M)-B$ of $M^*$.
\end{lemma}

The Tutte polynomial, as defined in equation~(\ref{tdef}), can
alternatively be expressed as follows:
\begin{equation}\label{tactive} t(M;x,y)= \sum_{B\in
\mathcal{B}(M)} x^{i(B)}y^{e(B)}.\end{equation} In particular,
although $i(B)$ and $e(B)$, for a particular basis $B$, depend on
the order $<$, the multiset of pairs $\bigl(i(B),e(B)\bigr)$, as
$B$ ranges over the bases of $M$, does not depend on the order.
Thus, the coefficient of $x^iy^j$ in $t(M;x,y)$ is the number of
bases of $M$ with internal activity $i$ and external activity $j$.

The crux of understanding the Tutte polynomial of a lattice path
matroid is describing internal and external activities of bases in
terms of the associated lattice paths; this is what we turn to
now.  Recall that if the bounding lattice paths $P$ and $Q$ go
from $(0,0)$ to $(m,r)$, then the lattice path matroid $M[P,Q]$
has ground set $[m+r]$; the elements in $[m+r]$ represent the
first step, the second step, and so on. We use the natural linear
order on $[m+r]$, that is, $1<2<\cdots<m+r$. We start with a lemma
that is an immediate corollary of Theorem~\ref{thm:bases}.

\begin{lemma}\label{iini}
Assume that $\{b_1,b_2,\ldots,b_r\}$ is a basis of a lattice path
matroid with $b_1<b_2<\cdots<b_r$.  Then $b_i$ is in the set $N_i$
of potential $i$-th North steps.
\end{lemma}

The following theorem describes externally active elements for
bases of lattice path matroids.

\begin{thm}
Assume that $B = \{b_1,b_2,\ldots,b_r\}$ is a basis of a lattice
path matroid $M[P,Q]$ with $b_1<b_2<\cdots<b_r$.  Assume that $x$
is not in $B$; say $b_i<x<b_{i+1}$. There is a $j$ with $j\leq i$
and with $(B - b_j)\cup x$ a basis of $M[P,Q]$ if and only if $x$
is in $N_i$. Equivalently, $x$ is externally active in $B$ if and
only if the $x$-th step of the lattice path that corresponds to
$B$ is an East step of the lower bounding path $P$.
\end{thm}

\begin{proof}
If $x$ is in $N_i$, then clearly $(B - b_i)\cup x$ is a
transversal of the set system $N_1,N_2,\ldots,N_r$ and so is a
basis of $M[P,Q]$. Conversely, if $(B - b_j)\cup x$ is a basis of
$M[P,Q]$ for some $j$ with $j\leq i$, then since $x$ is the $i$-th
element in this basis, Lemma~\ref{iini} implies that $x$ is in
$N_i$. The equivalent formulation of external activity follows
immediately by interpreting the first assertion in terms of
lattice paths.
\end{proof}

By the last theorem, Lemma~\ref{dualact}, and the lattice path
interpretation of matroid duality, we get the following result.

\begin{thm}\label{active}
Let $B$ be a basis of the lattice path matroid $M[P,Q]$ and let
$P(B)$ be the lattice path associated with $B$. Then $i(B)$ is the
number of times $P(B)$ meets the upper path $Q$ in a North step
and $e(B)$ is the number of times $P(B)$ meets the lower path $P$
in an East step.
\end{thm}

\begin{figure}
\begin{center}
\epsfxsize 2.25truein \epsffile{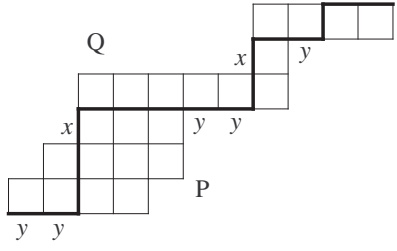}
\end{center}
\caption{The lattice path corresponding to a basis with internal
activity $2$ and external activity $5$, which contributes $x^2y^5$
to the Tutte polynomial.}\label{actfig}
\end{figure}

Theorem~\ref{active} is illustrated in Figure~\ref{actfig}. It is
worth noting the following simpler formulation in the case of
$k$-Catalan matroids.

\begin{cor}\label{kcatact}
Let $B$ be a basis of a $k$-Catalan matroid and let $P(B)$ be the
associated lattice path. Then $i(B)$ is the number of times $P(B)$
returns to the line $y=x/k$ and $e(B)$ is $j$ where $(j,0)$ is the
last point on the $x$-axis in $P(B)$.
\end{cor}

This lattice path interpretation of basis activities is one of the
keys for obtaining the following generating function for the
sequence of Tutte polynomials of the $k$-Catalan matroids.

\begin{thm}\label{thm:tuttegf}
Let $$C=C(z) = \sum_{n\geq 0} \frac{1}{kn+1}\binom{(k+1)n}{n}
z^n$$ be the generating function for the $k$-Catalan numbers. The
generating function for the Tutte polynomials of the $k$-Catalan
matroids is
\begin{equation}\label{tutteform}
 \sum_{n\geq 0} t(M_n^k;x,y)z^n=
1+\left(\frac{xzy^k}{1-z \sum_{l=1}^k y^l C^{k-l+1}}\right)
\frac{1}{1-xzC^{k}}.
\end{equation}
\end{thm}

\begin{proof}
From our lattice path interpretation of bases of the $k$-Catalan
matroid $M_n^k$, we are concerned with lattice paths that
\begin{itemize}
\item[(i)] go from $(0,0)$ to $(kn,n)$ and
\item[(ii)]  do not go above the line $y=x/k$.
\end{itemize}
We consider two special types of such lattice paths. Let $d_{jn}$
be the number of such lattice paths that, in addition, have the
following two properties:
\begin{itemize}
\item[(iii)] the last point of the path that is on the $x$-axis is the
           point $(j,0)$ and
\item[(iv)] the path returns to the line $y=x/k$ exactly once.
\end{itemize}
By property (iv), we have $d_{j0} = 0$ for all $j$. Let $D =
D(y,z)=\sum_{n,j>0} d_{jn}y^jz^n$. Let $e_{in}$ be the number of
lattice paths that satisfy properties (i)--(ii) and the following
property:
\begin{itemize}
\item[(iii$'$)] the path returns to the line $y=x/k$ exactly $i$ times.
\end{itemize}
Here the term $e_{00}$ is $1$.  Let $E = E(x,z)=\sum_{n,i\geq 0}
e_{in}x^iz^n$. By the lattice path interpretation of bases and
activities, we have
\begin{equation}\label{product}
\sum_{n\geq 0} t(M_n^k;x,y)z^n= 1 + x\,D(y,z)\,E(x,z).
\end{equation}
Equation~(\ref{tutteform}) follows immediately from
equation~(\ref{product}) and the following two equations, the
justifications of which are the focus of the rest of the proof.
\begin{equation}\label{E}
E(x,z)=\frac{1}{1-xz C(z)^k}
\end{equation}
\begin{equation}\label{D}
D(y,z)=\frac{z y^k}{1-z \sum_{l=1}^{k} y^l C(z)^{k-l+1}}
\end{equation}
To prove equations~(\ref{E}) and (\ref{D}), it will be convenient
to use the notation $l_s$ for the line $y=(x-s)/k$.

Equation~(\ref{E}) is immediate once we prove that the generating
function $\sum e_{1n} z^n$ for the number of paths that return
exactly once to the line $y=x/k$ is $zC(z)^k$.  To see this,
consider the following decomposition of such a path $P^*$. By
considering the last point of $P^*$ that is on the line $l_1$, and
then the last point of $P^*$ that is on $l_2$, and so on up to
$l_k$, it follows that the path $P^*$ can be decomposed uniquely
as a sequence $P^*=E P_1 E P_2\cdots E P_{k-1} E P_k N$, where
$P_i$ is a path beginning and ending on the line $l_i$ and never
going above this line.

We turn to equation~(\ref{D}).  Let $P^*$ be a path that returns
to the line $y=x/k$ exactly once.  If $P^*$ consists of $kn$ East
steps followed by $n$ North steps, then $P^*$ contributes $z^n
y^{kn}$ to $D(y,z)$; all such paths contribute $\sum_{i\geq 1}
z^iy^{ki}$, that is  $z y^k/(1-z y^k)$, to $D(y,z)$.  Assume path
$P^*$ is not of this type.  Let $i$ be the minimum value of $s$
such that $P^*$ intersects $l_s$ in a point neither on the line
$y=0$ nor on $x=kn$. Let $t$ be $\lceil \frac{i}{k} \rceil$. Since
$P^*$ contains the point $(kn,n-t)$, it follows that $P^*$ can be
decomposed uniquely as follows: $$P^*=E^i P'P_i E P_{i+1} E \cdots
E P_{t k } N^t,$$ where $P'$ is a non-empty path that begins and
ends on the line $l_i$ and that returns only once to this line,
and $P_s$ is a path that begins and ends on the line $l_s$ and
does not go above this line.  There are $kt-i+1$ paths among $P_i,
P_{i+1}, \ldots, P_{t k }$ and such paths are enumerated by
$C(z)$.  If $i$ is $kt$, then the path $P_i E P_{i+1} E \cdots E
P_{t k }$ reduces to $P_i$. In this case if the path $P_i$ were
trivial, then the path $P^*$ would intersect the line $l_i$ only
in the lines $y=0$ and $x=kn$, which contradicts the choice of
$i$. Therefore, when $i$ is $kt$ we have to guarantee that $P_i$
is nontrivial. Hence, we get
\begin{align}
 D = & \, \frac{zy^k}{1-zy^k} +
    \sum_{\substack{i\,:\,i\geq 1, \\ kt-i+1\ne 1}} y^i\, D\, z^t
C^{kt-i+1}
  + \sum_{\substack{i\,:\,i\geq 1, \\ kt-i+1 = 1}} y^i\, D\, z^t (C-1)
\notag \\
 = & \, \frac{zy^k}{1-zy^k} +
    D \sum_{i\geq 1} y^i z^t C^{kt-i+1}
  - D \sum_{\substack{i\,:\,i\geq 1, \\ kt-i+1 = 1}} y^i z^t . \notag
\end{align}
Since $kt-i+1$ is $1$ if and only if $i$ is $kt$, the last term
above is $D\,z y^k/(1-z y^k)$. To simplify the rest, note that
$\sum_{i\geq 1} y^i z^t C^{kt-i+1}$ is $$\begin{matrix} \quad\quad
y z\, C^k & + & \quad  y^2 z\, C^{k-1} &  + & \cdots & + &
\,\,\,\,y^k z\, C \,\, +  \\ y^{k+1} z^2 C^k & + & y^{k+2} z^2
C^{k-1} & +  &\cdots & + & y^{2k} z^2 C \,\,  +  \\ y^{2k+1} z^3
C^k & + & y^{2k+2} z^3 C^{k-1} & +  &\cdots  & + & y^{3k} z^3 C
\,\,   + \\ \quad\quad\vdots &  &\quad\quad\vdots  & & \cdots& &
\quad\quad \vdots\,\, \notag
\end{matrix}$$
which, by adding the columns, gives $$\frac{z}{1-z
y^k}\sum_{l=1}^k y^l  C^{k-l+1}.$$ Thus, $$D = \frac{zy^k}{1-zy^k}
+ \frac{D\,z}{1-z y^k}\sum_{l=1}^k  y^l C^{k-l+1} - \frac{D\,z
y^k}{1-z y^k}. $$ Solving for $D$ gives equation~(\ref{D}),
thereby completing the proof of the theorem.
\end{proof}

By extracting the coefficients of the expression in
Theorem~\ref{thm:tuttegf} we find a formula for the coefficients
of the Tutte polynomial of a $k$-Catalan matroid. To  write this
formula more compactly, let us denote by $S(m,s,k)$ the number of
solutions to the equation$$ l_1+\cdots +l_s=m$$ such that $1\leq
l_i\leq k$ for all $i$ with $1\leq i\leq s$. Set $S(0,0,k)=1$. It
will be useful to note that $$S(m,s,1)=\left\{\begin{array}{ll} 1,
& \mathrm{if }\  m=s; \\ 0, &\mathrm{otherwise.}
\end{array}\right. $$ An elementary inclusion-exclusion argument gives
 $$S(m,s,k)= \sum_{i=0}^s (-1)^i \binom{s}{i}\binom{m-ki-1}{s-1}.$$

\begin{thm}\label{thm:tuttecoefs}
The coefficient of $x^iy^j$ in the Tutte polynomial $t(M^k_n;x,y)$
of the $k$-Catalan matroid $M^k_n$ is $$\sum_{s=0}^{m} S(m,s,k)
\binom{(k+1)(n-1)-i-m}{n-s-i-1}\frac{s(k+1)-m+k(i-1)}{n-s-i},$$
where $m=j-k$. Equivalently, this is the number of lattice paths
that
\begin{itemize}
\item[(i)] go from $(0,0)$ to $(kn,n)$,

\item[(ii)] use steps $(1,0)$ and $(0,1)$,
\item[(iii)] do not go above the line $y = x/k$,
\item[(iv)] have as their last point on the $x$-axis the point $(j,0)$,
and
\item[(v)] return to the line $y = x/k$ exactly $i$ times.
\end{itemize}
\end{thm}

\begin{proof}
We need to extract the coefficient of $x^iy^jz^n$ in
equation~(\ref{tutteform}). We start by extracting the coefficient
of $y^{j-k}$ in $$\frac{1}{1-z \sum_{l=1}^k y^l
C^{k-l+1}}=\sum_{s\geq 0} \left(z \sum_{l=1}^k y^l
C^{k-l+1}\right)^s.$$ Let $m=j-k$. From basic algebra the
coefficient of $y^m$ in $(z \sum_{l=1}^k y^l C^{k-l+1})^s$ is
$z^sS(m,s,k)C^{s(k+1)-m}$. From this it follows that the
coefficient of $x^iy^j$ in the right-hand side of
equation~(\ref{tutteform}) is $$z^{i}C^{k(i-1)}\left(\sum_{s=0}^m
z^sS(m,s,k)C^{s(k+1)-m} \right). $$ To conclude the proof, we have
to extract the coefficient of $z^n$ in the above expression; this
is done using Lemma~\ref{lem:coeffcat}.
\end{proof}

It is an open problem to obtain explicit expressions for the Tutte
polynomials of the matroids $M_n^{k,l}$ for values of $k$ and $l$
not covered by the previous theorem, namely $k>1$ and $l>1$. The
first unsolved case is $k=l=2$. The sequence
$1,6,53,554,6362,77580,\dots$ that gives the number of bases of
$M_n^{2,2}$ also arises in the enumeration of certain types of
planar trees, and in that context Lou Shapiro gave a nice
expression for the corresponding generating function (see entry
A066357 in~\cite{eis}). This sequence also appears in~\cite{tbp};
indeed, as we show in Section~\ref{sec:tennis}, there is a simple
connection between the number of bases in certain lattice path
matroids and the problem considered in~\cite{tbp}.

We single out a corollary of Theorem~\ref{thm:tuttecoefs} that
shows a very rare property possessed by the Tutte polynomials of
the Catalan matroids $M_n$.

\begin{cor}
For $n>1$, the Tutte polynomial of the Catalan matroid $M_n$ is
$$\sum_{i,j>0}
 \frac{i+j-2}{n-1}\binom{2n-i-j-1}{n-i-j+1}x^iy^j.$$ In
particular, the coefficient of $x^i y^j$ in the Tutte polynomial
$t(M_n;x,y)$ of the Catalan matroid $M_n$ depends only on $n$ and
the sum $i+j$.
\end{cor}

We close this section with some simple observations.  A well-known
corollary of Lemma~\ref{dualact} is that the Tutte polynomial of a
matroid and its dual are related by the following equation:
$$t(M^*;x,y) = t(M;y,x).$$ From this and Theorem~\ref{thm:dualkl}
we get the following corollary.

\begin{cor}
The Tutte polynomial of the $(k,l)$- and $(l,k)$-Catalan matroids
are related as follows: $$t(M^{k,l}_n;x,y) = t(M^{l,k}_n;y,x).$$
Thus, the Tutte polynomial of the $(k,k)$-Catalan matroid
$M^{k,k}_n$, and in particular the Catalan matroid $M_n$, is a
symmetric function in $x$ and $y$.
\end{cor}

\section{Computing the Tutte polynomial of lattice path
matroids}\label{sec:computing}

There is no known polynomial-time algorithm for computing the
Tutte polynomial of an arbitrary matroid, or even its evaluations
at certain points in the plane~\cite{ckcc}. There are many
evaluations of the Tutte polynomial that are particularly
significant; for instance, it follows from
equation~(\ref{tactive}) that $t(M;1,1)$ is the number of bases of
$M$. Since the bases of a lattice path matroid correspond to paths
that stay in a given region and the number of such paths is given
by a determinant (see Theorem~$1$ in Section~$2$.$2$
of~\cite{moh}), the number of bases of a lattice path matroid can
be computed in polynomial time. It turns out that other
evaluations like $t(M;1,0)$ and $t(M;0,1)$ can also be expressed
as determinants. This led us to suspect  that the Tutte polynomial
of a lattice path matroid could be computed in polynomial time. In
this section, we show that this is indeed the case: we give such a
polynomial-time algorithm. Also, we give a second technique for
computing the Tutte polynomial in the case of generalized Catalan
matroids (this second technique, although more limited in scope,
is particularly simple to implement using standard mathematical
software). The results in this section stand in striking contrast
to those in~\cite{cpv}, where it is shown that for fixed $x$ and
$y$ with $(x-1)(y-1)\ne 1$, the problem of computing $t(M;x,y)$
for a transversal matroid $M$ is $\#$P-complete.

\begin{figure}
\begin{center}
\epsfxsize 3.0truein \epsffile{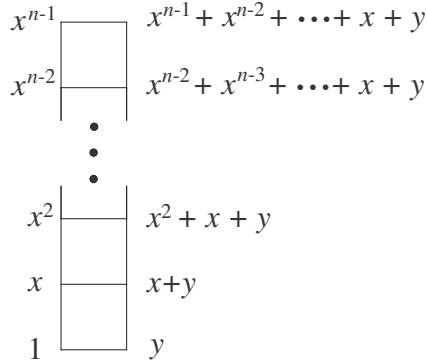}
\end{center}
\caption{The recursive computation of the Tutte polynomial of an
$n$-element circuit via lattice path statistics.}\label{cir}
\end{figure}

By Theorem~\ref{active}, for a lattice path matroid $M=M[P,Q]$,
the Tutte polynomial $t(M;x,y)$ is the generating function
$$\sum_{B\in\mathcal{B}(M)} x^{i(B)}y^{e(B)}$$ where $i(B)$ is the
number of North steps that the lattice path $P(B)$ corresponding
to $B$ shares with the upper bounding path $Q$ and $e(B)$ is the
number of East steps that $P(B)$ shares with the lower bounding
path $P$. Any lattice path can be viewed as a sequence of shorter
lattice paths.  This perspective gives the following algorithm for
computing the Tutte polynomial of the lattice path matroid
$M=M[P,Q]$ where $P$ and $Q$ go from $(0,0)$ to $(m,r)$. With each
lattice point $(i,j)$ in the region $R$ bounded by $P$ and $Q$,
associate the polynomial $$f(i,j) = \sum_{P'} x^{i(P')}y^{e(P')}$$
where the sum ranges over the lattice paths $P'$ that go from
$(0,0)$ to $(i,j)$ and stay in the region $R$, and where, as for
$t(M;x,y)$, the exponent $i(P')$ is the number of North steps that
$P'$ shares with $Q$ and $e(P')$ is the number of East steps that
$P'$ shares with $P$.  In particular, $f(m,r) = t(M;x,y)$.  Note
that for a point $(i,j)$ in $R$ other than $(0,0)$, at least one
of $(i-1,j)$ or $(i,j-1)$ is in $R$; furthermore, only $(i-1,j)$
is in $R$ if and only if the step from $(i-1,j)$ to $(i,j)$ is an
East step of $P$, and, similarly, only $(i,j-1)$ is in $R$ if and
only if the step from $(i,j-1)$ to $(i,j)$ is a North step of $Q$.
The following rules for computing $f(i,j)$ are evident from these
observations and the definition of $f(i,j)$.
\begin{itemize}
\item[(a)] $f(0,0) = 1$.
\item[(b)] If the lattice points $(i,j)$, $(i-1,j)$ and $(i,j-1)$
are all in the region $R$, then $f(i,j)= f(i-1,j)+f(i,j-1)$.
\item[(c)] If the lattice points $(i,j)$ and $(i-1,j)$ are in $R$
but $(i,j-1)$ is not in $R$, then $f(i,j)= y\,f(i-1,j)$.
\item[(d)] If the lattice points $(i,j)$ and $(i,j-1)$ are in $R$
but $(i-1,j)$ is not in $R$, then $f(i,j)= x\, f(i,j-1)$.
\end{itemize}
This algorithm is illustrated in Figure~\ref{cir} where we apply
it to compute the Tutte polynomial of an $n$-element circuit.

If $x$ and $y$ are set to $1$, the algorithm above reduces to a
well-known technique for counting lattice paths. This is
consistent with the general theory of Tutte polynomials; as noted
above, $t(M;1,1)$ is the number of bases of $M$.

The recurrence above requires at most $(r+1)(m+1)$ steps to
compute the Tutte polynomial of a lattice path matroid whose
bounding paths go from $(0,0)$ to $(m,r)$. Thus, we have the
following corollary.

\begin{cor}
The Tutte polynomial of a lattice path matroid can be computed in
polynomial time.
\end{cor}

We remark that the recurrence expressed in (a)--(d) above is
essentially the deletion-contraction rule for Tutte polynomials,
along with the corresponding rules for loops and isthmuses (see,
e.g.,~\cite{bo} for this perspective on the Tutte polynomial).
This follows by considering the lattice path interpretations of
deletion and contraction, which are given in~\cite{lpm2}. We also
remark that while the deletion-contraction rule for computing
$t(M;x,y)$ generally gives rise to a binary tree with $2^{|E(M)|}$
leaves, for lattice path matroids there are relatively few
isomorphism types for minors, and the geometry of lattice paths
automatically collects minors of the same isomorphism type. To
make this more specific, let $R$ be the region bounded by the
lattice paths $P$ and $Q$ of the lattice path matroid $M=M[P,Q]$.
As can be seen from the description of minors in~\cite{lpm2}, each
minor whose ground set is an initial segment $[k]$ of $[m+r]$ can
be viewed as having as bases the lattice paths in $R$ from $(0,0)$
to some specific point in $R$ of the form $(i,k-i)$. It follows
that the number of possible minors of $M=M[P,Q]$ that arise when
computing $t(M;x,y)$, rather than being exponential, is bounded
above by $(r+1)(m+1)$.

By Theorem~\ref{thm:catrec}, generalized Catalan matroids are
formed from the empty matroid by iterating the operations of
taking free extensions and direct sums with the uniform matroid
$U_{1,1}$.  The following rule is well-known and easy to check:
for any matroid $M$,
\begin{equation}\label{ti}
t(M\oplus U_{1,1};x,y) = x \, t(M;x,y).
\end{equation}
For free extensions, we have the following result, which is easy
to prove using formula~(\ref{tdef}) and the rank function of the
free extension. (This formula is equivalent to the expression for
the Tutte polynomial of a free extension given in Proposition 4.2
of~\cite{bry1}.)

\begin{thm}
The Tutte polynomial of the free extension $M+e$ of $M$ is given
by the formula
\begin{equation}\label{tfe}
t(M+e;x,y) = \frac{x}{x-1}\,t(M;x,y) +
\left(y+\frac{x}{x-1}\right)t(M;1,y).
\end{equation}
\end{thm}

Formulas~(\ref{ti}) and~(\ref{tfe}) can be used, for instance, to
compute Tutte polynomials of $(k,l)$-Catalan matroids very
quickly. It is through such computations that we were lead, for
instance, to Theorem~\ref{thm:betaud}.

\section{The broken circuit complex and the characteristic
polynomial}\label{sec:char}

In this section we study two related objects for lattice path
matroids, the broken circuit complex and the characteristic
polynomial. The second of these is an invariant of the matroid but
the first depends on a linear ordering of the elements. We show
that under the natural ordering of the elements, the broken
circuit complex of any loopless lattice path matroid has a
property that is not shared by the broken circuit complexes of
arbitrary matroids, namely, the broken circuit complex of a
lattice path matroid is the independence complex of another
matroid, indeed, of a lattice path matroid. Our study of the
characteristic polynomial is more specialized; we focus on the
characteristic polynomial $\chi(\widehat{M^k_n};\lambda)$ of the
matroid $\widehat{M^k_n}$ obtained from the $k$-Catalan matroid
$M^k_n$ by omitting the loops.  Our results on the broken circuit
complex lead to a lattice path interpretation of each coefficient
of $\chi(\widehat{M^k_n};\lambda)$ from which we obtain a formula
for these coefficients.  We start by outlining the necessary
background on broken circuit complexes; for an extensive account,
see~\cite{shell}.

Given a matroid $M$ and a linear order $<$ on the ground set
$E(M)$, a {\em broken circuit} of the resulting ordered matroid is
a set of the form $C-x$ where $C$ is a circuit of $M$ and $x$ is
the least element of $C$ relative to the linear ordering. A subset
of  $E(M)$ is an {\em nbc-set} if it contains no broken circuit.
Clearly subsets of nbc-sets are nbc-sets. Thus, $E(M)$ and the
collection of nbc-sets of $M$ form a simplicial complex, the {\em
broken circuit complex of $M$ relative to $<$}, which is denoted
$BC_<(M)$.  Different orderings of $E(M)$ can produce
nonisomorphic broken circuit complexes (see,
e.g.,~\cite[Example~7.4.4]{shell}). The facets of $BC_<(M)$ are
the {\em nbc-basis}, that is, the basis of $M$ that are nbc-sets.
The following characterization of nbc-bases is well-known and easy
to prove.

\begin{lemma}\label{lem:refext}
The nbc-basis of $M$ are the bases of $M$ of external activity
zero.
\end{lemma}

Note that nbc-sets contain no circuits and so are independent.
Thus, the broken circuit complex $BC_<(M)$ of $M$ is contained in
the {\em independence complex of} $M$, that is, the complex with
ground set $E(M)$ in which the faces are the independent sets of
$M$.  Note also that if $M$ has loops, then the empty set is a
broken circuit, so $M$ has no nbc-sets.  Thus, throughout this
section we consider only matroids with no loops.

As in Section~\ref{sec:tutte}, we use the natural ordering on the
points of lattice path matroids. The examples in~\cite{shell} show
that the broken circuit complex need not be the independence
complex of another matroid. In contrast, Theorem~\ref{thm:bccx}
shows that the broken circuit complex of a lattice path matroid
without loops is the independence complex of another lattice path
matroid.

\begin{thm}\label{thm:bccx}
With the natural order, the broken circuit complex of a lattice
path matroid $M[P,Q]$ with no loops is the independence complex of
the lattice path matroid $M[P',Q]$ where $NP=P'N$.
\end{thm}

\begin{proof}
Since a subset of $E(M[P,Q])$ is an nbc-set of $M[P,Q]$ if and
only if it is contained in an nbc-basis of $M[P,Q]$, it suffices
to show that the nbc-bases of $M[P,Q]$ are precisely the bases of
$M[P',Q]$.  By Lemma~\ref{lem:refext} and Theorem~\ref{active},
the nbc-bases of $M[P,Q]$ correspond to the lattice paths in the
region bounded by $P$ and $Q$ that share no East step with $P$.
Thus, the nbc-bases of $M[P,Q]$ correspond to the lattice paths in
the region bounded by $P'$ and $Q$ where the East steps of $P'$
occur exactly one unit above those of $P$.  This condition on $P'$
is captured by the equality $NP=P'N$.
\end{proof}

\begin{cor}\label{cor:charnbcsets}
Let $M[Q]$ be a generalized Catalan matroid with no loops. A
subset $X$ of the ground set of $M[Q]$ is an nbc-set if and only
if $X\cup 1$ is independent in $M[Q]$. In particular, all
independent sets of $M[Q]$ that contain $1$ are nbc-sets and the
nbc-basis of $M$ are exactly the bases of $M[Q]$ that contain $1$.
\end{cor}

We now turn to the characteristic polynomial, which plays an
important role in many enumeration problems in matroid theory
(see~\cite{rota,zas}) and which can be defined in a variety of
ways.  As mentioned above, the isomorphism type of the broken
circuit complex of a matroid $M$ can depend on the ordering of the
points.  However, it can be shown that the number of nbc-sets of
each size is an invariant of the matroid; these numbers are the
coefficients of the characteristic polynomial. Specifically, the
\emph{characteristic polynomial} $\chi(M;\lambda)$ of a matroid
$M$ is
\begin{equation}\label{nbccharpoly}
\chi(M;\lambda)= \sum_{i=0}^{r(M)}
 (-1)^i \, \mathbf{nbc}(M;i)\, {\lambda}^{r(M)-i},
\end{equation}
where $\mathbf{nbc}(M;i)$ is the number of nbc-sets of size $i$.
Thus, $(-1)^{r(M)}\chi(M;-\lambda)$ is the face enumerator of the
broken circuit complex of $M$.  (Equation~(\ref{nbccharpoly})
applies even if the matroid $M$ has loops, in which case
$\chi(M;\lambda)$ is $0$.) Alternatively, $\chi(M;\lambda)$ can be
expressed in terms of the Tutte polynomial as follows:
 $$\chi(M;\lambda) = (-1)^{r(M)}
t(M;1-\lambda,0)=\sum_{A\subseteq E(M)} (-1)^{|A|}
\lambda^{r(M)-r(A)}.$$
 The characteristic polynomial can also be expressed in the
following way in terms of the M\"obius function of the lattice of
flats:
 $$\chi(M;\lambda) =
\sum\limits_{\substack{\text{ flats }  F \\
  \text{ of } M }} \mu(\emptyset,F) \lambda^{r(M)-r(F)}.$$
 (See, e.g.,~\cite[Theorem 7.4.6]{shell},
for details.) In particular,  the absolute value of the constant
term of $\chi(M;\lambda)$ is both the number of nbc-bases of $M$
and the absolute value of the M\"obius function $\mu(M)$. This and
Theorem~\ref{thm:bccx} give the following corollary.

\begin{cor}\label{cor:mu}
The absolute value of the M\"obius function $\mu(M[P,Q])$ of a
loopless lattice path matroid is the number of bases of the
lattice path matroid $M[P',Q]$ where $NP=P'N$, or, equivalently,
of the lattice path matroid $M[P^*,Q^*]$ where $P=P^*N$ and
$Q=NQ^*$.
\end{cor}

Our interest is in the characteristic polynomial of a specific
type of lattice path matroid. Recall that the elements
$1,2,\ldots,k$ are loops of the $k$-Catalan matroid $M^k_n$. Thus,
the characteristic polynomial of $M^k_n$ is zero. This motivates
considering the {\em loopless Catalan matroid} $\widehat{M_n}$,
which we define to be $M[(NE)^{n-1}N]$, and more generally the
{\em loopless $k$-Catalan matroid} $\widehat{M^k_n}$, which we
define to be $M[(NE^k)^{n-1}N]$. Thus, these matroids are formed
from almost the same bounding paths as those for the Catalan and
$k$-Catalan matroids except that the initial East steps that give
loops have been omitted.

We start with the following consequence of Corollary~\ref{cor:mu}.

\begin{cor}\label{mu}
The number of nbc-bases of $\widehat{M^k_n}$, that is,
$|\mu(\widehat{M^k_n})|$, is the $k$-Catalan number $C^k_{n-1}$.
In particular, $|\mu(\widehat{M_n})| = C_{n-1}$.
\end{cor}

\begin{proof}
From the second part of Corollary~\ref{cor:mu}, we have that
$|\mu(\widehat{M^k_n})|$ is the number of bases of $M^k_{n-1}$,
which is $C^k_{n-1}$.
\end{proof}

By combining Corollary~\ref{cor:indpaths} and
Corollary~\ref{cor:charnbcsets}, we get the following
characterization of the nbc-sets of size $i$ of $\widehat{M^k_n}$
in terms of lattice paths.

\begin{cor}\label{cor:nbcpaths}
Via the map $X\mapsto P(X)$, the nbc-sets of size $i$, for $0\leq
i \leq n$, in the loopless $k$-Catalan matroid $\widehat{M^k_n}$
correspond bijectively to the following two types of lattice
paths.
\begin{itemize}
\item[(i)] Lattice paths from $(0,1)$ to $\bigl((k+1)(n-1)-i+1,i\bigr)$
 that do not go above the line $y = \frac{1}{k}x+1$.
\item[(ii)] Lattice paths from $(0,1)$ to
 $\bigl((k+1)(n-1)-i,i+1\bigr)$ that do not go above the line
 $y = \frac{1}{k}x+1$.
\end{itemize}
\end{cor}

By using this characterization of nbc-sets we obtain the following
expression for each coefficient of the characteristic polynomial.

\begin{thm}\label{thm:coefscar}
The absolute value of the coefficient of ${\lambda}^{n-i}$ in the
characteristic polynomial of the loopless $k$-Catalan matroid
$\widehat{M_n^k}$ is given by the formula
 $$\mathbf{nbc}(\widehat{M_n^k};i) = \left\{\begin{array}{ll}1, &
\mathrm{if}\ i=0;\\ \frac{\displaystyle
(k+1)(n-i-1)+2}{\displaystyle (k+1)(n-1)+2}{\displaystyle
\binom{(k+1)(n-1)+2}{i}},& \mathrm{if}\ 1\leq i \leq n-1;\\
C^k_{n-1}, & \mathrm{if}\ i=n.
\end{array}\right.$$
\end{thm}

\begin{proof}
Since $\widehat{M^k_n}$ is loopless, the empty set is an nbc-set;
from this the case $i=0$ follows.  The case $i=n$ has been treated
in Corollary~\ref{mu}. For $i$ with $1\leq i\leq n-1$, we have to
count the number of paths as in Corollary~\ref{cor:nbcpaths}. This
is equivalent to counting the following:
\begin{itemize}
\item[(i)] lattice paths from $(0,0)$ to $\bigl((k+1)(n-1)-i+1,i-1\bigr)$
 that do not go above the line $y = x/k$, and
\item[(ii)] lattice paths from $(0,0)$ to
 $\bigl((k+1)(n-1)-i,i\bigr)$ that do not go above the line
 $y = x/k$.
\end{itemize}
Observe that the sum of the number of paths described in items (i)
and (ii) is the number of paths from $(0,0)$ to
$\bigl((k+1)(n-1)-i+1,i\bigr)$ that do not go above the line
$y=x/k$. The formula follows then from Lemma~\ref{ballot}.
\end{proof}

From the formula in Theorem~\ref{thm:coefscar} and appropriate
manipulation, we see that the linear term in the characteristic
polynomial of $\widehat{M_n}$ is also a Catalan number. Note that,
however, for the loopless $k$-Catalan matroid the linear term of
the characteristic polynomial is not the corresponding $k$-Catalan
number.

\begin{cor}\label{linearterm}
The linear term in $\chi (\widehat{M_n},\lambda)$ is $C_n$.
\end{cor}

\section{The $\beta$ invariant}\label{sec:beta}

The $\beta$ invariant $\beta(M)$ of a matroid $M$, which was
introduced by Crapo, can be defined in several ways;
see~\cite[Section 3]{zas} for a variety of perspectives on the
$\beta$~invariant, as well as its applications to connectivity and
series-parallel networks.  We use the following definition.  It
can be shown that for any matroid $M$, the coefficients of $x$ and
$y$ in the Tutte polynomial $t(M;x,y)$ are the same; this
coefficient is $\beta(M)$. Since loops are externally active with
respect to every basis, no basis of a matroid $M$ with loops will
have external activity zero, so $\beta(M)$ is zero; dually, if $M$
has isthmuses, then $\beta(M)$ is zero. Therefore, in this section
we focus on matroids with neither loops nor isthmuses.

Let $N^{k,k}_{n}$ be the generalized Catalan matroid whose upper
path is $Q=(N^kE^k)^n$.  It is clear from the lattice path
presentation that $N^{k,k}_{n}$ is formed from the $(k,k)$-Catalan
matroid $M^{k,k}_{n+1}$ by deleting the $k$ loops and the $k$
isthmuses. The main result of this section is that
$\beta(N^{k,k}_{n})$ is $k$~times the Catalan number $C_{kn-1}$.
This result was suggested by looking at examples of Tutte
polynomials of lattice path matroids, but it can be formulated
entirely in terms of lattice paths, which is the perspective we
use in the proof. Indeed, the result is most striking when viewed
in terms of lattice paths.

The $\beta$~invariant of $N^{k,k}_{n}$ is the number of bases with
internal activity one and external activity zero; let $B$ be such
a basis and let $P(B)$ be its associated lattice path. By
Theorem~\ref{active}, the first step of $P(B)$ is $N$, the second
is $E$, and $P(B)$ does not contain any other North step in $Q$.
It is easy to see that such lattice paths $P(B)$ are in $1$-$1$
correspondence with the paths from $(0,0)$ to $(kn-1,kn-1)$ that
do not go above the path $N^{k-1}(E^kN^k)^{n-1}E^{k-1}$. Recall
that the number of paths from $(0,0)$ to $(kn-1,kn-1)$ that do not
go above the line $y=x$ is $C_{kn-1}$.  In this section we show
that the number of paths that do not go above the path
$N^{k-1}(E^kN^k)^{n-1}E^{k-1}$ is $k$ times $C_{kn-1}$. We start
with the case $k=1$.

\begin{thm}\label{thm:betacat}
The $\beta$ invariant of $N^{1,1}_n$ is $C_{n-1}$.
\end{thm}

\begin{proof}
By the discussion above, $\beta(N^{1,1}_n)$ is the number of paths
from $(0,0)$ to $(n-1,n-1)$ that do not go above the path
$(EN)^{n-1}$, which is $C_{n-1}$.
\end{proof}

From here on, we consider only paths that use steps $U$ and $D$.
From the discussion above and the correspondence between the
alphabets, we get the following lemma.

\begin{lemma}\label{lem:betaud}
The $\beta$~invariant of the matroid $N^{k,k}_n$ is the number of
paths that
\begin{itemize}
\item[(i)] go from $(0,0)$ to $\bigl(2(nk-1),0\bigr)$,
\item[(ii)] use steps $U$ and $D$, and
\item[(iii)] never go below the path
$D^{k-1}(U^kD^k)^{n-1}U^{k-1}$.
\end{itemize}
\end{lemma}

\begin{figure}
\epsfxsize 12truecm \epsffile{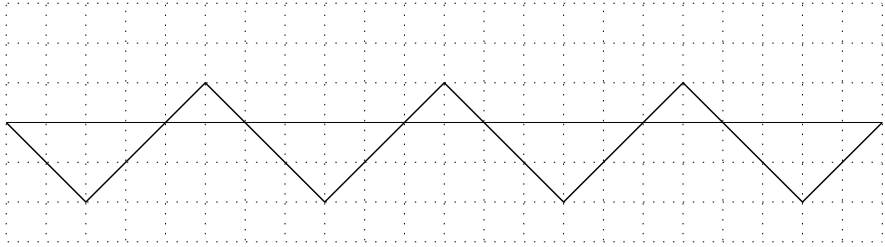} \caption{The path
$D^{k-1}(U^kD^k)^{n-1}U^{k-1}$ for $k=3$ and
$n=4$.}\label{fig:betaud}
\end{figure}

A path of the form $D^{k-1}(U^kD^k)^{n-1}U^{k-1}$ is depicted in
Figure~\ref{fig:betaud}. The next theorem is the main result of
this section.

\begin{thm}\label{thm:betaud}
The number of paths that go from $(0,0)$ to
$\bigl(2(nk-1),0\bigr)$, use steps $U$ and $D$, and do not go
below the path $D^{k-1}(U^kD^k)^{n-1}U^{k-1}$ is $k\,C_{nk-1}$.
\end{thm}

Before proving the theorem, we mention that  if we change the
bounding path to  $(D^kU^k)^n$, the elegance and brevity of the
result seem to disappear; currently there is no known comparably
simple answer. Indeed, the path $(D^kU^k)^n$ is connected with an
open problem in enumeration that is discussed in the next section.
The following corollary is an immediate consequence of
Lemma~\ref{lem:betaud} and Theorem~\ref{thm:betaud}.

\begin{cor}
The $\beta$~invariant of the matroid $N^{k,k}_n$ is $k\,C_{nk-1}$.
\end{cor}

\noindent \textit{Proof of Theorem~\ref{thm:betaud}.} Let us
denote the path $D^{k-1}(U^kD^k)^{n-1}U^{k-1}$ by $B$. In this
proof we consider paths from $(0,0)$ to $(2kn-1,-1)$ using steps
$U$ and $D$. When we say that one such path does not go below a
given border, we refer to the path with the last step removed.
Hence, a Dyck path is a path from $(0,0)$ to $(2kn-1,-1)$ that
does not go below the line $y=0$ (except for the last step). A
\emph{cyclic permutation} of a path $s_1s_2\dots s_{l}$ is a path
$s_is_{i+1}\dots s_{l} s_1\dots s_{i-1}$ for some $i$ with $1\leq
i\leq l$. It is easy to show that all cyclic permutations of a
Dyck path from $(0,0)$ to $(2kn-1,-1)$ are different; note that
this does not hold if we consider Dyck paths ending in a point of
the form $(2l,0)$.

The proof is in the spirit of several results generically known as
the Cycle Lemma (see the notes at the end of Chapter~5
of~\cite{ec2}). One such result states that among the $2l+1$
possible cyclic permutations of a path from $(0,0)$ to
$(2l+1,-1)$, there is exactly one that is a Dyck path. Moreover,
the cyclic permutation that leads to a Dyck path is the one that
starts after the leftmost minimum of the path (see~\cite[Theorem
1.1]{stanton} for more details on this). To prove the theorem, we
show that for every Dyck path from $(0,0)$ to $(2kn-1,-1)$,
exactly $k$ of its cyclic permutations are paths that do not go
below $B$; conversely, every path that does not go below $B$ can
be obtained as one of these $k$ cyclic permutations of a Dyck
path. To describe these permutations we need to introduce some
terminology.

\begin{figure}
\epsfxsize 12truecm \epsffile{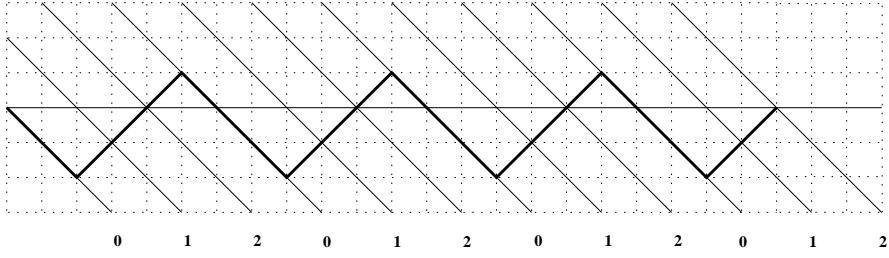} \caption{The partition
of the points $(x,y)$ for which $x+y$ is even
($k=3$).}\label{fig:classes}
\end{figure}

It is clear that if a lattice point $(x,y)$ is in a path that
begins at $(0,0)$ and uses steps $U$ and $D$, then  $x+y$ is even.
We partition the lattice points whose coordinates have an  even
sum into $k$ disjoint classes. The point $(x,y)$ is \emph{in
class} $c$ with $0\leq c\leq k-1$ if $(x+y)/2 \equiv c$ modulo
$k$. As can be seen in Figure~\ref{fig:classes}, each class
corresponds to an infinite family of parallel lines.

We say that a point $(x,y)$ has \emph{height} $y$. It is easy to
see that a point in class $c$ is not below the path $B$ if and
only if the height of the point is strictly greater than $c-k$.
Let $p=(x,y)$ be a point in a path that uses steps $U$ and $D$; we
say that $p$ is a \emph{down point} if $p$ is the end of a $D$
step. The \emph{cyclic permutation at} $p$ is the permutation that
starts in the step that has $p$ as the first point.

\begin{figure}
\epsfxsize 12truecm \epsffile{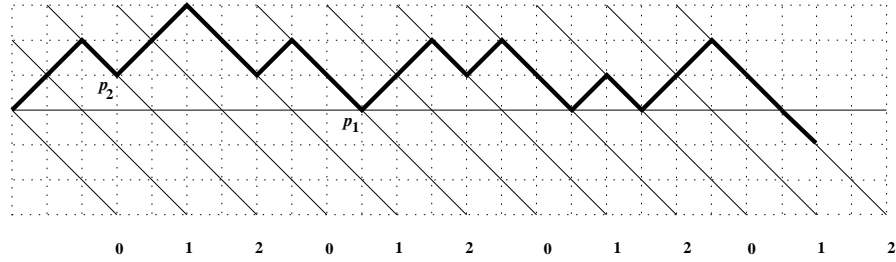} \caption{A Dyck path
and the points $p_i$ of the proof of Theorem~\ref{thm:betaud}.}
\label{fig:betadyck}
\end{figure}

Let $R$ be a Dyck path from $(0,0)$ to $(2kn-1,-1)$; clearly, $R$
does not go below $B$. The other $k-1$ cyclic permutations of $R$
that do not go below $B$ are given by the points
$p_1,\dots,p_{k-1}$ that we  define next. The point $p_1$ is the
first down point of $R$ that is in class $k-1$ and has height $0$.
The point $p_2$ is the first down point of $R$ that is in class
$k-2$ and has height $0$, if such a point exists; otherwise, take
the first down point in class $k-1$ and with height $1$. To find
the $i$-th point $p_i$, among all down points that are in class
$k-i+j$ and have height $j$ for $0\leq j\leq i-1$, take the ones
that have minimum $j$, and among those take as $p_i$ the one that
appears first in $R$. See Figure~\ref{fig:betadyck} for an
example.

We next show that the points $p_i$ exist. Since $R$ is a Dyck
path, the point $\bigl(2(kn-1),0\bigr)$ is always in $R$.
Moreover, it is a down point and belongs to class $k-1$; hence
there is at least one down point in $R$ in class $k-1$ with height
$0$. In general, we prove that if for $j$ with $j<i-1$ there is no
down point in class $k-i+j$ with height $j$, then there exists a
down point in class $k-1$ with height $i-1$, and thus we take as
$p_i$ the first such point that appears in $R$. Assume that $R$
contains no point in class $k-i+j$ with height $j$ for $j$ with
$j<i-1$. Since $R$ is a Dyck path all points at height 0 are down,
so the point $\bigl(2(kn-i),0\bigr)$ is not in the path $R$. The
point $\bigl(2(kn-i)+1,1\bigr)$ is in class $k-i+1$ and has height
$1$, so by assumption it cannot be down. Then, if it is in $R$,
the previous step must be $U$, but that forces the point
$\bigl(2(kn-i),0\bigr)$ to be in $R$, which is a contradiction.
Hence $\bigl(2(kn-i)+1,1\bigr)$ is not in $R$. In the same way one
proves that the points of the form $\bigl(2(kn-i)+j,j\bigr)$ are
not in $R$ for $j$ with $0\leq j\leq i-2$. Note that this implies
that the path $R$ goes above all these points. Now consider the
point $\bigl(2(kn-i)+i-1,i-1\bigr)$, which is in class $k-1$. If
this point is not in $R$, then the point in $R$ with first
coordinate equal to $2(kn-i)+i-1$ would have height at least
$i+1$; however, from such a point it is impossible to reach the
point $(2(kn-1),0)$, which is always in the path. Therefore, the
point $\bigl(2(kn-i)+i-1,i-1\bigr)$ is in $R$, and since the point
$\bigl(2(kn-i)+i-2,i-2\bigr)$ is not, it has to be a down point.
Hence $R$ contains a down point in class $k-1$ with height $i-1$,
and the existence of $p_i$ is proved.

Now we have to check that $\pi_i (R)$, the cyclic permutation of
$R$ at $p_i$, is a path that does not go below the path $B$. We
split $\pi_i(R)$ into two subpaths $R_1$ and $R_2$ such that
$R=R_1R_2$ and $\pi_i (R)=R_2R_1$. We prove that there is no point
in either part $R_1$ or $R_2$ of $\pi_i (R)$ below the path $B$.
Assume that $p_i$ belongs to class $k-i+j$ and has height $j$ for
some $j$ with $j\leq i-1$.

Suppose first there is a point in the subpath $R_1$ that goes
below $B$ and let $q$ be the first such point; this point is a
down point and if it is in class $c$, then its height is $c-k$.
Let us move the point $q$ to $R$, that is, let the point $q_R$ be
the point of $R$ that goes to $q$ after the cyclic permutation at
$p_i$. It is easy to check that the point $q_R$ has height
$j+1+c-k$ and belongs to class $c+j-i+1$ modulo $k$. Since $R$ is
a Dyck path we have $j+1+c-k\geq 0$; from this and the inequality
$j\leq i-1$ it follows that the class of $q_R$ is indeed
$c+j-i+1$. Since $c<k$, we have that $j+1+c-k\leq j$. This
together with the fact that the point $q_R$ comes before $p_i$ in
$R$ contradict the choice of $p_i$.

Similarly, suppose there is a point in the subpath $R_2$ of
$\pi_i(R)$ that goes below $B$ and let $q'$ be the first such
point. As before, the point $q'$ is down and if it is in class
$c'$, then its height is $c'-k$. Let $q'_R$ be the point of $R$
that is mapped to $q'$ by the cyclic permutation at $p_i$. The
point $q'_R$ has height $j+c'-k$; thus since $R$ is a Dyck path,
$j+c'-k\geq 0$. The class of $q'_R$ is $k-i+j+c'$ modulo $k$. By
combining the inequalities $j\leq i-1$, $c'<k$, and $j+c'-k\geq
0$, we get that class of $q_R$ is $k-i+(j+c'-k)$. Since
$j+c'-k<j$, the point $q_R$ contradicts the choice of $p_i$. This
finishes the proof that the cyclic permutation at $p_i$ is a path
that does not go below $B$.

We now have that every Dyck path from $(0,0)$ to $(2kn-1,-1)$
gives rise to $k$ paths that do not go below $B$, including the
Dyck path itself. As noted above, all cyclic permutations of a
Dyck path are different and for every path only one cyclic
permutation is a Dyck path. Since there are $C_{kn-1}$ Dyck paths,
we have that the number of paths as described in the statement of
the theorem is at least $k\,C_{kn-1}$.

To complete the proof of the equality, we have to show that every
path that does not go below $B$ is either a Dyck path $R$ or one
of the $k-1$ cyclic permutations of a Dyck path $R$ at one of the
points $p_1,p_2,\dots,p_{k-1}$ defined above. Let $S$ be a path
from $(0,0)$ to $(2kn-1,-1)$ that does not go below $B$ and that
is not a Dyck path; let $q_0$ be its first point and $q_S$ its
leftmost minimum. The cyclic permutation at $q_S$ is a Dyck path
$S'$. Let $q'_0$ be the image of the point $q_0$ in $S'$. If the
point $q_S$ is in class $c$ and has height $h$ in $S$, then the
point $q'_0$ in $S'$ is a down point that belongs to class $k-1-c$
and has height $-h-1$; also $h>c-k$. We have to show that $q'_0$
is one of the points $p_1,\dots,p_{k-1}$ with respect to the Dyck
path $S'$. Since by definition the point $p_i$ is in class $k-i+j$
and has height $j$, it follows that $q'_0$ should be the point
$p_{c-h}$ with $j=-h-1$ (note that $c-h$ is a valid index since
$c-k<h$). The result will follow if we show that no down point in
$S'$ is in class $k-c+h+j$ and has height $j$ for $0\leq j<-h-1$,
and that any down point in class $k-c-1$ with height $-h-1$ comes
after $q'_0$ in $S'$. It is easy to show that if there were a
point satisfying either condition, then its height in the path $S$
would exceed the class minus $k$, and hence the point would be
below the path $B$, which is a contradiction. \hfill $\Box$

\section{Connections with the tennis ball problem}\label{sec:tennis}

The following problem is of current interest in enumerative
combinatorics; only a very limited number of cases have been
settled (see~\cite{tbp}).

\medskip

\noindent{\sc The $(k+l,l)$ tennis ball problem.}  Let
$b_1,b_2,\ldots,b_{(k+l)n}$ be a sequence of distinct balls. At
stage $1$, balls $b_1,b_2,\ldots,b_{k+l}$ are put in bin $A$ and
then $l$ balls are moved from bin $A$ to bin $B$. At stage $i$,
balls $b_{(i-1)(k+l)+1},b_{(i-1)(k+l)+2},\ldots,b_{i(k+l)}$ are
put in bin $A$ and  then some set of $l$ balls from bin $A$ are
moved to bin $B$. (In particular, balls that remain in bin $A$
after stage $i-1$ can go in bin $B$ at stage $i$.) How many
different sets of $ln$ balls can be in bin $B$ after $n$
iterations?

\medskip

We show that the answer is the number of bases of the
$(k,l)$-Catalan matroid $M^{k,l}_{n+1}$.

It is well known that free extensions of transversal matroids are
also transversal; we use the presentations of free extensions
given in the following lemma.

\begin{lemma}\label{repfree}
Assume that $M$ is a transversal matroid of rank $r$ with
presentation $(A_j:j\in K)$ where $|K|=r$.  Then the free
extension $M+e$ is also transversal and the set system $(A_j\cup
e:j\in K)$ is a presentation of $M+e$.
\end{lemma}

\begin{proof}
The partial transversals $X$ of $(A_j\cup e:j\in K)$ with
$e\not\in X$ are precisely the partial transversals of $(A_j:j\in
K)$. Also, for any partial transversal $X$ of $(A_j:j\in K)$ with
$|X|<r$, the set $X\cup e$ is a partial transversal of $(A_j\cup
e:j\in K)$.
\end{proof}

There are many ways to add a set $I=\{f_1,f_2,\ldots,f_u\}$ of
isthmuses to a transversal matroid with presentation $(A_j:j\in
K)$; for instance, $(A_j:j\in K)$ together with
$\{f_1\},\{f_2\},\ldots,\{f_u\}$ is such a presentation. The
presentation of interest for us is the union of the multiset
$(A_j\cup I :j\in K)$ with $u$ copies of $I$.

\begin{figure}
\begin{center}
\epsfxsize 5.5truein \epsffile{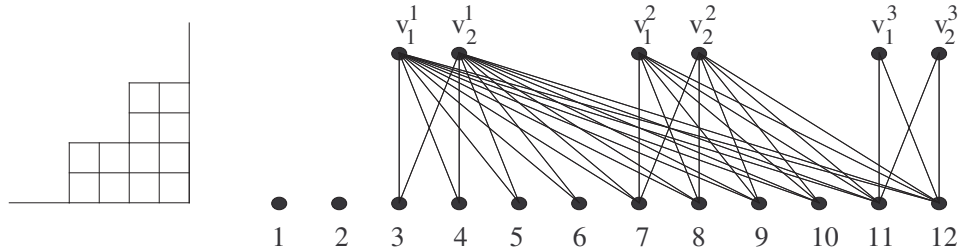}
\end{center}
\caption{The lattice path presentation of $M^{2,2}_3$ and the
graph $G^{2,2}_3$.}\label{graph}
\end{figure}

By Theorem~\ref{thm:catrec}, the $(k,l)$-Catalan matroid
$M^{k,l}_{n+1}$ can be constructed from the empty matroid by
taking $k$ free extensions, then adding $l$ isthmuses, then taking
$k$ free extensions, then adding $l$ isthmuses, etc., for a total
of $n+1$ iterations.  With this in mind, as well as the
presentations of free extensions and additions of isthmuses just
discussed, consider the following bipartite graph $G^{k,l}_{n+1}$.
One set of the bipartition of the vertex set is $[(k+l)(n+1)]$,
the ground set of $M^{k,l}_{n+1}$; let $v^j_h$, with $1\leq j \leq
n+1$ and $1 \leq h \leq l$, be the remaining vertices.  Vertices
$(k+l)i+\kappa$, with $1\leq \kappa \leq k$, are adjacent to all
$v^j_h$ with $1\leq j \leq i$ and $1 \leq h \leq l$; vertices
$(k+l)i+\eta$, with $k+1\leq \eta \leq k+l$, are adjacent to all
$v^j_h$ with $1\leq j \leq i+1$ and $1 \leq h \leq l$. The graph
$G^{2,2}_3$ is illustrated in Figure~\ref{graph}.

It follows from the descriptions of presentations of free
extensions and extensions by isthmuses that the bases of
$M^{k,l}_{n+1}$ are precisely the sets of vertices in
$[(k+l)(n+1)]$ of maximal size that can be matched in
$G^{k,l}_{n+1}$. Note that $M^{k,l}_{n+1}$ has as many bases as
the matroid obtained by deleting the first $k$ elements (which are
loops) and the last $l$ elements (which are isthmuses); let
$\widehat{G}^{k,l}_{n+1}$ denote graph obtained from
$G^{k,l}_{n+1}$ by deleting these vertices. The graph
$\widehat{G}^{k,l}_{n+1}$ can be used to model $n$ iterations of
the $(k+l,l)$-tennis ball problem: after relabelling vertices,
those adjacent to $v^{n}_1, v^{n}_2,\ldots,v^{n}_l$ can be viewed
as the balls that could be selected to go in bin $B$ on first
iteration; those adjacent to
$v^{n-1}_1,v^{n-1}_2,\ldots,v^{n-1}_l$ can be viewed as the
candidates to go in bin $B$ on the second iteration, and so on.
Furthermore, maximal-sized sets of vertices that can be matched in
this graph are precisely the sets of balls that can be in bin $B$
at the end of $n$ iterations. Thus, the answer to the $(k+l,l)$
tennis ball problem, with $n$ iterations, is the number of bases
of the $(k,l)$-Catalan matroid $M^{k,l}_{n+1}$.

\begin{center}
\textsc{Acknowledgements}
\end{center}

The  authors thank Lou Shapiro for interesting and useful
discussions about the Catalan numbers and the objects they count,
and for bringing the tennis ball problem to our attention. The
first author also thanks Bill Schmitt for further discussions
about Catalan numbers.

\end{document}